\documentclass[12pt,a4paper]{article}
\usepackage{amsmath}
\usepackage{amssymb}
\usepackage{theorem}
\usepackage{graphicx}

\numberwithin{equation}{section}    

\theoremstyle{plain}            

\newtheorem{te}{Theorem}[section]

\newtheorem{pr}[te]{Proposition}

\theoremstyle{definition}       

\newtheorem{re}[te]{Remark}

\usepackage[dvips]{color}
\usepackage{eucal}
\usepackage{amssymb,amsmath,amsfonts}
\usepackage{indentfirst}
\usepackage{natbib} 
\date{}
\begin{document}

\title{\textbf{The Dynamics of Rabinovich system  }}

 \author{ \textbf{Oana Chi\c{s} and Mircea Puta$^{\dag}$}}
\maketitle

Abstract: The paper presents some dynamical aspects of Rabinovich
type, with distributed delay and with fractional derivatives.
\\* 2000 AMS Mathematics Subject Classification: 37K10,
26A33, 58A05, 58A40, 53D17
\\* Keywords: Rabinovich system, Poisson
representation, metriplectic structure, distributed delay,
fractional derivatives, Caputo fractional derivatives, Caputo
integral operator.
\section{Introduction}
\par As an important application of chaotic dynamical systems, chaos-based
secure communication and cryptography attracted continuous
interest over the last decade. It studies methods of controlling
deterministic systems with chaotic behavior. Moreover, it is easy
to notice the possibility of substantial variation of the
characteristics of chaotic systems by relatively small variations
of their parameters and external actions. A method of transmitting
information using chaotic signal was proposed by A.S.Pikovsky and
M.I. Rabinovich \cite{Pi} using the differential system
$$\left\{%
\begin{array}{ll}
    \dot{x_{1}}=-\nu_{1}x_{1}+hx_{2}+x_{2}x_{3} &  \\
    \dot{x_{2}}=hx_{1}-\nu_{2}x_{2}-x_{1}x_{3} &  \\
    \dot{x_{3}}=-\nu_{3}x_{3}+x_{1}x_{2}.&  \\
\end{array}%
\right.$$
\par In this paper we will consider the Rabinovich system:
\begin{equation}\label{1.1}
\left\{%
\begin{array}{ll}
    \dot{x_{1}}=x_{2}x_{3} &  \\
    \dot{x_{2}}=-x_{1}x_{3} &  \\
    \dot{x_{3}}=x_{1}x_{2}&  \\
\end{array},%
\right.
\end{equation}
and we will analyze some global properties, the local study of
stationary points, compatible Poisson structures and corresponding
tri-Hamiltonian systems are also discussed.
\par A Hamiltonian equation is called tri-Hamiltonian if
     it admits two Hamiltonian representations with compatible
     Poisson structures $$\frac{dx}{dt}=J\nabla H=\tilde{J}\nabla \tilde{H}=\bar{J}\nabla
     \bar{H},$$ where $J, \, \tilde{J}$ and $\bar{J}$ are three
     Hamiltonian matrices (of the form $[\{x_{i},x_{j}\}]$ and
     $\{\cdot,\cdot\}$ is the Poisson structure) and they are also compatible.

\par The Hamiltonian formulation is important in mathematics, physics and also in other branches of natural science. From
another point of view, we considered the analysis of the revised
dynamical system, the analysis of the dynamical system with
distributive delay variables and the analysis of the fractional
dynamical system.
\par Revised dynamical system, with distributive delay,
associated to system (\ref{1.1}), allow the description of new
crypting methods.

\section{The analysis of classical Rabinovich differential
equations}
\subsection{\textbf{Geometrical properties
of the system (\ref{1.1})}}
\par In this subsection we will present some dynamical and
geometrical properties, from geometrical mechanical point of view,
(\cite{Ch},\cite{Pu2}).
\begin{pr}\label{0} The dynamics
 (\ref{1.1}) have the following Hamilton-Poisson
realizations:
  \par(i) $(\mathbb{R}^{3}, P^{i}, h_{i}), \quad i=1, 2, 3$ where
\begin{eqnarray*} &P^{1}&=\begin{bmatrix}
  0 & x_{3} &-x_{2} \\
  -x_{3} & 0 &  0\\
  x_{2} & 0 & 0 \\
\end{bmatrix},\quad
h_{1}(x_{1}, x_{2}, x_{3})= \frac{1}{2}(x_{1}^{2}+x_{2}^{2});\\
 &P^{2}&=\begin{bmatrix}
  0 & 0 &\frac{1}{2}x_{2} \\
  0 & 0 &  -\frac{1}{2}x_{1}\\
  -\frac{1}{2}x_{2} & \frac{1}{2}x_{1} & 0 \\
\end{bmatrix}, \quad
h_{2}(x_{1}, x_{2}, x_{3})= x_{2}^{2}+x_{3}^{2}; \\
 &P^{3}&=\begin{bmatrix}
  0 & 0 &-\frac{1}{2}x_{2} \\
  0 & 0 &  \frac{1}{2}x_{1}\\
  \frac{1}{2}x_{2} & -\frac{1}{2}x_{1} & 0 \\
\end{bmatrix},\quad
h_{3}(x_{1}, x_{2}, x_{3})= x_{1}^{2}-x_{3}^{2}.\quad
\end{eqnarray*}
   \par(ii) $(\mathbb{R}^{3}, P_{123}^{\alpha \beta \gamma }, h_{\alpha}),$ where
$$P_{123}^{\alpha \beta \gamma }=\alpha P_{1}+\beta P_{2}+\gamma
P_{3}=\begin{bmatrix}
  0 & \alpha x_{3} & -(\alpha-\frac{\beta}{2}+\frac{\gamma}{2})x_{2} \\
  -\alpha x_{3} & 0 & -(\frac{\beta}{2}-\frac{\gamma}{2})x_{1}\\
  (\alpha-\frac{\beta}{2}+\frac{\gamma}{2})x_{2} & (\frac{\beta}{2}-\frac{\gamma}{2})x_{1} & 0 \\
\end{bmatrix},\\$$ and
$h_{\alpha}(x_{1},x_{2},x_{3})=\frac{1}{2\alpha}(x_{1}^{2}+x_{2}^{2}),$
for each $\alpha, \beta, \gamma \in \mathbb{R}, \, \alpha\neq 0.$
$\Box$\end{pr}
\par From direct computations, using the algebraic technique of  Bermejo and
Fairen \cite{Be}, we get the following results.
\begin{pr}
There exists only one functionally independent Casimir of our
Poisson configurations:
\begin{itemize}
    \item $(\mathbb{R}^{3}, P_{1})$  given by:
$c_{1}(x_{1}, x_{2}, x_{3})=\frac{1}{2}(x_{2}^{2}+x_{3}^{2});$
    \item $(\mathbb{R}^{3}, P_{2})$  given by:
$c_{2}(x_{1}, x_{2}, x_{3})=x_{1}^{2}+x_{2}^{2};$
    \item $(\mathbb{R}^{3}, P_{3})$  given by:
$c_{3}(x_{1}, x_{2}, x_{3})=x_{1}^{2}+x_{2}^{2};$
    \item $(\mathbb{R}^{3},
P_{123}^{\alpha\beta\gamma})$ given by:
$c_{\alpha\beta\gamma}(x_{1}, x_{2},
x_{3})=-\frac{1}{\alpha}(\frac{\beta}{2}-\frac{\gamma}{2})x_{1}^{2}+\frac{1}{\alpha}(\alpha-
\frac{\beta}{2}+\frac{\gamma}{2})x_{2}^{2}+ x_{3}^{2}. $
\end{itemize}$\Box$
\end{pr}

\subsection {\textbf{Stability problem}}
\par From the analysis of the stationary points, using \cite{Pu1}, we
get the following statements.
\begin{pr} The stationary points $e_{1}^{m}(m,0,0),$ $e_{2}^{m}(0,m,0)$ and
$e_{3}^{m}(0,0,m)$, $m\in \mathbb{R}$ have the following behavior:
\par (i) $e_{1}^{m}(m,0,0),$ $m\in \mathbb{R}$ are spectrally stable;
\par (ii) $e_{2}^{m}(0,m,0),$ $m\in \mathbb{R}$ are unstable; \par (iii)
$e_{3}^{m}(0,0,m),$ $m\in \mathbb{R}$ are spectrally stable.$\Box$
\end{pr}

\begin{pr}
The stationary points $e_{1}^{m}(m,0,0)$ and $e_{3}^{m}(0,0,m)$
are nonlinear stable.$\Box$
\end{pr}
\par Now we will point out periodic orbits and
heteroclinic orbits associated to Rabinovich system. When studying
the existence of periodic solution, we will use \cite{Bi1}.
\begin{pr} The dynamics reduced to the coadjoint orbit
$(x_{1})^{2}+(x_{2})^{2}=m^{2}$ has near the stationary points
$e_{1}^{m}(m,0,0)$ $m\in \mathbb{R}^{*}$ at least one periodic
solution whose period is close to $\frac{\pi}{\mid m\mid}.$$\Box$
\end{pr}
\begin{pr} There exists four heteroclinic orbits between the
stationary points $e_{2}^{m}(0, m, 0)$ and $e_{2}^{-m}(0, -m, 0)$,
$m\in \mathbb{R}, \, m\neq 0$ given by:
\begin{equation}
\left\{%
\begin{array}{ll}
    x_{1}(t)=\pm m\, sech(mt) &  \\
    x_{2}(t)=\pm m\,tgh(mt) &  \\
    x_{3}(t)=\pm m\,sech(mt).&  \\
\end{array}%
\right.
\end{equation}These orbits belong to the planes $x_{3}=
\pm x_{1}$.$\Box$
\end{pr}

\section{The metriplectic structure associated to \\ Rabinovich system}
\par A Leibniz structure on a smooth manifold M is defined by a
tensor field P of type (2,0), \cite{Al3}. The tensor field P and a
smooth function h on M, called a Hamiltonian function, define a
vector field $X_{h}$ which generates a differential system, called
Leibniz system. If P is skew-symmetric then we have an almost
simplectic structure and if P is symmetric then we have an almost
metric structure.
\par Let P be a skew-symmetric tensor field in $\mathbb{R}^{3}$ of
type (2,0), g a 2-symmetric tensor field and $h\in C^{\infty}
(\mathbb{R}^{3}).$ If P is a Poisson tensor field and g is a
nondegenerate tensor field, then $(\mathbb{R}^{3},P,g)$ is called
a metriplectic manifold of the first kind (\cite{Bi2}, \cite{Fi},
\cite{Or}). The differential system is given by
\begin{equation}\label{3.1}
\dot{x}_{i}=\sum_{j=1}^{3}P_{ij}\frac{\partial h}{\partial
x_{j}}+\sum_{j=1}^{3}g_{ij}\frac{\partial h}{\partial x_{j}},
\quad i=1,2,3.
\end{equation}
\par The $g_{ij}$ is compatible with h, and is given by:

\begin{equation}\label{3.2}
g_{ii}=-\sum_{k=1,k\neq i}^{3}\frac{\partial h}{\partial
x_{k}}\frac{\partial h}{\partial x_{k}}, \quad
g_{ij}=\frac{\partial h}{\partial x_{i}}\frac{\partial h}{\partial
x_{j}}, \quad i,j=1,2,3.
\end{equation}

\par If P is a (almost) Poisson differential system on
$\mathbb{R}^{3}$ with Hamiltonian function $h_{1}$ and a Casimir
function $h_{2}$, there exists a tensor field g such that
$(\mathbb{R}^{3},P,g)$ is a metriplectic manifold of second kind.
The differential system associated with it is given by:
\begin{equation}\label{3.3}
\dot{x}_{i}=\sum_{j=1}^{3}P_{ij}\frac{\partial h_{1}}{\partial
x_{j}}+\sum_{j=1}^{3}g_{ij}\frac{\partial h_{2}}{\partial x_{j}},
\quad i=1,2,3,
\end{equation} where

\begin{equation}\label{3.4}
g_{ii}=-\sum_{k=1,k\neq i}^{3}\frac{\partial h_{1}}{\partial
x_{k}}\frac{\partial h_{2}}{\partial x_{k}}, \quad
g_{ij}=\frac{\partial h_{1}}{\partial x_{i}}\frac{\partial
h_{2}}{\partial x_{j}}, \quad i,j=1,2,3.
\end{equation}
\par Let $(\mathbb{R}^{3},P^{\alpha}), \, \alpha=1,2,3$ realizations of Rabinovich system of differential
equations, with Hamiltonian functions $h_{\alpha}, \,
\alpha=1,2,3$ and Casimir functions $c_{\alpha}, \, \alpha=1,2,3,$
where:
\begin{equation}\label{3.5}
 P^{1}=\begin{bmatrix}
  0 & x_{3} &-x_{2} \\
  -x_{3} & 0 &  0\\
  x_{2} & 0 & 0 \\
\end{bmatrix},\\ \quad h_{1}=\frac{1}{2}(x_{1}^{2}+x_{2}^{2}), \quad
c_{1}=\frac{1}{2}(x_{2}^{2}+x_{3}^{2});
\end{equation}

\begin{equation}\label{3.6}
 P^{2}=\begin{bmatrix}
  0 & 0 &\frac{1}{2}x_{2} \\
  0 & 0 &  -\frac{1}{2}x_{1}\\
  -\frac{1}{2}x_{2} & \frac{1}{2}x_{1} & 0 \\
\end{bmatrix},\\ \quad h_{2}=\frac{1}{2}(x_{1}^{2}+x_{2}^{2}), \quad
c_{2}=\frac{1}{2}(x_{2}^{2}+x_{3}^{2});
\end{equation}

\begin{equation}\label{3.7}
 P^{3}=\begin{bmatrix}
 0 & 0 & -\frac{1}{2}x_{2} \\
  0 & 0 & \frac{1}{2}x_{1}\\
  \frac{1}{2}x_{2} & -\frac{1}{2}x_{1} & 0 \\
\end{bmatrix},\\ \quad h_{3}=\frac{1}{2}(x_{1}^{2}+x_{2}^{2}), \quad
c_{3}=\frac{1}{2}(x_{2}^{2}+x_{3}^{2}).
\end{equation}

\par Using (\ref{3.1}), (\ref{3.2}), (\ref{3.3}) and (\ref{3.4}) we get
the following results.

\begin{pr}
\par (a)The metriplectic realization of the first kind of (\ref{3.5}) is
given by $(\mathbb{R}^{3},P_{1},g_{1})$ where:
$$g_{11}^{1}=-(x_{2})^{2}, \quad g_{22}^{1}=-(x_{1})^{2}, \quad g_{33}^{1}=0;$$
$$g_{12}^{1}=x_{1}x_{2}, \, g_{21}^{1}=x_{1}x_{2}, \, g_{13}^{1}=g_{31}^{1}=0, \, g_{23}^{1}=g_{32}^{1}=0.$$

\par (b) The associated differential system is given by:
\begin{equation}\label{3.8}
\left\{%
\begin{array}{ll}
    \dot{x_{1}}=x_{2}x_{3}+x_{1}x_{2}(x_{1}-x_{2}) &  \\
    \dot{x_{2}}=-x_{1}x_{3}+(x_{1})^{2} &  \\
    \dot{x_{3}}=x_{1}x_{2}.&  \\
\end{array}%
\right.
\end{equation}

\par(c)The differential system (\ref{3.8}) has the following  stationary
points: $$e_{1}^{m}(m,0,0), \quad e_{2}^{m}(0,m,0), \quad
e_{3}^{m}(0,0,m).$$

\par(d) The matrix of the linear part of the system (\ref{3.8}) in $e_{1}^{m}(m,0,0)$,
$e_{2}^{m}(0,m,0)$, resp in $e_{3}^{m}(0,0,m)$ is given by:$$
A_{1}=\begin{bmatrix}
 0 & m^{2} & 0 \\
 0 & 0 & m^{2}+m\\
 0 & m & 0 \\
\end{bmatrix}, \quad A_{2}=\begin{bmatrix}
 -m^{2} & 0 & m \\
 0 & 0 & 0\\
 m & 0 & 0 \\
\end{bmatrix}, \quad resp \quad A_{3}=\begin{bmatrix}
 0 & m & 0 \\
 -m & 0 & 0\\
 0 & 0 & 0 \\
\end{bmatrix}.
$$

\par(e) The characteristic equation of $A_{1}$ for (\ref{3.9}) in
$e_{1}^{m}(m,0,0)$ is:
$$\lambda(-\lambda^{2}+m^{2}(m+1))=0 $$and so, we have two cases:
\par (i) if $m>-1$, then $e_{1}^{m}(m,0,0)$ are unstable; \par (ii) if $m<-1$,
then we have a limit cycle.

\par(f) The characteristic equation of $A_{2}$ for (\ref{3.8}) in
$e_{2}^{m}(0,m,0)$ is:
$$-\lambda(\lambda^{2}+m^{2}\lambda-m^{2})=0$$ and so, it can be
easily seen that $e_{2}^{2}(0,m,0)$ are unstable.

\par(g) The characteristic equation of $A_{3}$ for (\ref{3.8}) in
$e_{3}^{3}(0,0,m)$ is:
$$\lambda(\lambda^{2}+m^{2})=0.$$

\par(h) In a neighborhood of $e_{3}^{m}(0,0,m),\, m>0$ there exists a limit
cycle.$\Box$
\end{pr}

\begin{pr}
\par (a) The metriplectic realization of the second kind of (\ref{3.5}) is
given by $(\mathbb{R}^{3},P_{1},g_{1})$ where:
$$g_{11}^{1}=-(x_{2})^{2}, \quad g_{22}^{1}=0, \quad g_{33}^{1}=0;$$
$$g_{12}^{1}=x_{1}x_{2}, \, g_{21}^{1}=0, \, g_{13}^{1}=x_{1}x_{3}, \, g_{31}^{1}=0, \, g_{23}^{1}=x_{2}x_{3}, \, g_{32}^{1}=0.$$

\par (b) The associated differential system is given by:
\begin{equation}\label{10}
\left\{%
\begin{array}{ll}
    \dot{x_{1}}=x_{2}x_{3}+x_{1}((x_{2})^{2}+(x_{3})^{2}) &  \\
    \dot{x_{2}}=-x_{1}x_{3}+x_{2}x_{3} &  \\
    \dot{x_{3}}=x_{1}x_{2}.&  \\
\end{array}%
\right.
\end{equation}

\par(c) The differential system (\ref{10}) has the following  stationary
points: $$e_{1}^{m}(m,0,0), \quad e_{2}^{m}(0,m,0), \quad
e_{3}^{m}(0,0,m).$$

\par(d) The matrix of the linear part of the system (\ref{3.8}) in $e_{1}^{m}(m,0,0)$,
$e_{2}^{m}(0,m,0)$, resp $e_{3}^{m}(0,0,m)$ is given by:$$
A_{1}=\begin{bmatrix}
 0 & 0 & 0 \\
 0 & 0 & -m\\
 0 & m & 0 \\
\end{bmatrix}, \quad A_{2}=\begin{bmatrix}
 m^{2} & 0 & m \\
 0 & 0 & m\\
 m & 0 & 0 \\
\end{bmatrix}, \quad resp \quad A_{3}=\begin{bmatrix}
 m^{2} & m & 0 \\
 -m & m & 0\\
 0 & 0 & 0 \\
\end{bmatrix}.
$$

\par(e) The characteristic equation of $A_{1}$ for (\ref{10}) in
$e_{1}^{m}(m,0,0)$ is:
$$\lambda(\lambda^{2}+m^{2})=0.$$

\par(f)The characteristic equation of $A_{2}$ for (\ref{10}) in
$e_{2}^{m}(0,m,0)$ is:
$$\lambda(\lambda^{2}-\lambda m^{2}-m^{2})=0 $$ and so, it can be
easily seen that $e^{m}_{2}(0,m,0)$ are unstable.

\par(g)The characteristic equation of $A_{3}$ for (\ref{10}) in
$e_{3}^{m}(0,0,m)$ is:
$$\lambda(\lambda^{2}-\lambda(m+m^{2})+m^{2})=0.$$ $\Box$
\end{pr}

\begin{re}
\par In an analogous way we can discuss the metriplectic
realization of first kind of (\ref{3.6}) and (\ref{3.7}).
\end{re}

\section{The differential systems with distributed delay}

\par Let us consider the product $\mathbb{R}^{3} \times \mathbb{R}^{3}=\{(\tilde{x},x)\mid \tilde{x}\in \mathbb{R}^{3}, x \in
\mathbb{R}^{3}\}$ and the canonical projections $\pi_{i} :
\mathbb{R}^{3}\times \mathbb{R}^{3}\longrightarrow \mathbb{R}^{3},
\quad i=1,2.$ A vector field $X\in \mathcal{X} (\mathbb{R}^{3}
\times \mathbb{R}^{3}),$ satisfying the condition
$X(\pi_{1}^{*}f)=0,$ for any $f\in
\mathbf{C}^{\infty}(\mathbb{R}^{3}),$ is given by:

\begin{equation}\label{4.1}
X(\tilde{x},x)=\sum_{i=1}^{n}X_{i}(\tilde{x},x)\frac{\partial}{\partial
x_{i}}.
\end{equation}
\par The differential system associated to X is given by:
\begin{equation}\label{4.2} \dot{x_{i}}(t)=X(\tilde{x},x), \quad
i=1,2,3.
\end{equation}

\par A differential system with distributed delay, see \cite{Al2} is a
differential system associated to a vector field $X\in \mathcal{X}
(\mathbb{R}^{3} \times \mathbb{R}^{3})$ for which
$X(\pi_{1}^{*}f)=0, \, \forall \, f\in
\mathbf{C}^{\infty}(\mathbb{R}^{3}),$ and it is given by
(\ref{4.1}) where $\tilde{x}(t)$ is:

\begin{equation}\label{4.3}
\tilde{x}(t)=\int_{0}^{\infty}k(s)x(t-s)ds
\end{equation}

k(s) is a distribution density. In the following we will consider
the following densities:

\begin{enumerate}
    \item uniform:
\begin{equation}\label{4.4}
k^{N}_{\tau}(s)=\left\{%
\begin{array}{ll}
    0, \quad 0\leq s \leq a &  \\
    \frac{1}{\tau}, \quad a\leq s \leq a+\tau &\\
    0, \quad s> a+\tau.&  \\
\end{array}%
\right.
\end{equation}
where $a>0$, $\tau>0$ are fixed numbers.
  \item exponential:
\begin{equation}\label{4.5}
k_{\alpha}(s)=\alpha e^{-\alpha s}, \quad \alpha>0;
\end{equation}
   \item Erlang:
   \begin{equation}\label{4.6}
k_{\alpha}(s)=\alpha^{2}s e^{-\alpha s}, \quad \alpha>0;
\end{equation}
    \item Dirac:
\begin{equation}\label{4.7}
k_{\alpha}(s)=\delta(s-\tau), \quad \tau>0;
\end{equation}
\end{enumerate}

\par The differential equations with distributed delay for
Rabinovich system are generated by an antisymmetric tensor field P
on $\mathbb{R}^{3}\times\mathbb{R}^{3}$ that  satisfies the
following relations:$$P(\pi_{1}^{*}f_{1}, \pi_{1}^{*}f_{2})=0,
\quad P(\pi_{2}^{*}f_{1}, \pi_{2}^{*}f_{2})=0$$ for all $f_{1},
f_{2}\in\mathbf{C}^{\infty}(\mathbb{R}^{3}).$

\par The differential
equation with distributed delay is given by:
\begin{equation}\label{4.9}
\dot{x}(t)=P(\tilde{x(t)},x(t))\nabla_{x}h(\tilde{x}(t),x(t)),
\end{equation}
where $\tilde{x}(t)=\int_{0}^{\infty}k(s)x(t-s)ds,$ and $h\in
\mathbf{C}^{\infty}(\mathbb{R}^{3}\times \mathbb{R}^{3}).$
\par Let
\begin{eqnarray*} P_{0}(x)&=&\begin{bmatrix}
  0 & x_{3} &-x_{2} \\
  -x_{3} & 0 &  0\\
  x_{2} & 0 & 0 \\
\end{bmatrix}, \, P_{1}(\tilde{x},x)=\begin{bmatrix}
  0 & x_{3} &-\tilde{x}_{2} \\
  -x_{3} & 0 &  0\\
  \tilde{x}_{2} & 0 & 0 \\
\end{bmatrix},\\
P_{2}(\tilde{x},x)&=&
\begin{bmatrix}
  0 & \tilde{x}_{3} &-x_{2} \\
  -\tilde{x}_{3} & 0 &  0\\
  x_{2} & 0 & 0 \\
\end{bmatrix},
P_{3}(\tilde{x},x)=\begin{bmatrix}
  0 & \tilde{x}_{3} &-\tilde{x}_{2} \\
  -\tilde{x}_{3} & 0 &  0\\
  \tilde{x}_{2} & 0 & 0 \\
\end{bmatrix}.\\
\end{eqnarray*}
\par We define
\begin{equation}\label{4.10}
P(\tilde{x},x)=\sum_{i=0}^{3}\varepsilon_{i}P_{i}, \quad with
\quad \varepsilon_{i}\geq 0, \quad
\sum_{i=0}^{3}\varepsilon_{i}=1.
\end{equation}
\par Let
\begin{eqnarray*}
h_{0}(\tilde{x},x)&=&\frac{1}{2}((x_{1})^{2}+(x_{2})^{2}),
h_{1}(\tilde{x},x)=\tilde{x}_{1}x_{1}+\frac{1}{2}(x_{2})^{2},\\
h_{2}(\tilde{x},x)&=&\frac{1}{2}(x_{1})^{2}+\tilde{x}_{2}x_{2},
h_{3}(\tilde{x},x)=\tilde{x}_{1}x_{1}+\tilde{x}_{2}x_{2}.\\
\end{eqnarray*}
\par We define
\begin{equation}\label{4.11}
h(\tilde{x},x)=\sum_{i=0}^{3}\delta_{i}h_{i}, \quad with \quad
\delta_{i}\geq 0, \quad \sum_{i=0}^{3}\delta_{i}=1.
\end{equation}
\par The Rabinovich differential equation with distributed delay
is given by (\ref{4.9}) with P and h given above by (\ref{4.10})
and (\ref{4.11}) with initial value $x(s)=\phi(s), \, s\in
(-\infty,0]$ where $\phi : (-\infty,0]\longrightarrow
\mathbb{R}^{3}, \, \phi\in \mathbf{C}^{\infty}(\mathbb{R}^{3}).$
\par In what follows we consider the functions $l\in
\mathbf{C}^{\infty}(\mathbb{R}^{3}\times\mathbb{R}^{3})$ given by:
\begin{eqnarray*}
l_{0}(\tilde{x},x)&=&\frac{1}{2}((x_{2})^{2}+(x_{3})^{2}),
l_{1}(\tilde{x},x)=\tilde{x}_{2}x_{2}+\frac{1}{2}(x_{3})^{2},\\
l_{2}(\tilde{x},x)&=&\frac{1}{2}(x_{2})^{2}+\tilde{x}_{3}x_{3},
l_{3}(\tilde{x},x)=\tilde{x}_{2}x_{2}+\tilde{x}_{3}x_{3}.\\
\end{eqnarray*}
\par We define
\begin{equation}\label{4.12}
l(\tilde{x},x)=\sum_{i=0}^{3}\varepsilon_{i}h_{i}, \quad with
\quad \varepsilon_{i}\geq 0, \quad
\sum_{i=0}^{3}\varepsilon_{i}=1.
\end{equation}
\begin{pr}
\begin{itemize}
    \item The function $l(\tilde{x},x)$  given by (\ref{4.12})
    satisfies the following relation:
\begin{equation}\label{4.13}
\nabla_{x}l(\tilde{x},x)P(\tilde{x},x)\nabla_{x}f(\tilde{x},x)=0,
\quad f\in
\mathbf{C}^{\infty}(\mathbb{R}^{3}\times\mathbb{R}^{3});
\end{equation}
    \item The revised differential equations with distributed
    delay satisfies the following relation:
\begin{equation}\label{4.14}
\dot{x}(t)=P(\tilde{x},x)\nabla_{x}h(\tilde{x},x)+g(\tilde{x},x)\nabla_{\tilde{x}}l(\tilde{x},x)
\end{equation} where $\tilde{x}(t)=\int_{0}^{\infty}k(s)x(t-s)ds$
and $g(\tilde{x},x)$ is a 2-tensor field given by:
\begin{eqnarray*}
g(\tilde{x},x)&=&(g_{ij}(\tilde{x},x)),\\
g_{ij}(\tilde{x},x)&=&\frac{\partial h(\tilde{x},x)}{\partial x_{i}}\frac{\partial h(\tilde{x},x)}{\partial x_{j}}, \, i\neq j\\
g_{ij}(\tilde{x},x)&=&-\sum_{k=i,k\neq i}\Big(\frac{\partial h(\tilde{x},x)}{\partial x_{k}}\Big )^{2}.\\
\end{eqnarray*}
\end{itemize}$\Box$
\end{pr}

\par The revised Rabinovich system with distributed delay has the
following form:
\begin{equation}\label{4.14}
\left\{%
\begin{array}{ll}
\dot{x_{1}}=(\alpha_{1}x_{3}+\alpha_{4}\tilde{x}_{3})(\beta_{1}x_{2}+\beta_{4}\tilde{x}_{2})+(\beta_{2}x_{1}+\beta_{3}\tilde{x}_{1})(\beta_{1}x_{2}+\beta_{4}\tilde{x}_{2})\alpha_{2}x_{3} &  \\
\dot{x_{2}}=-(\alpha_{1}x_{3}+\alpha_{4}\tilde{x}_{3})(\beta_{2}x_{1}+\beta_{3}\tilde{x}_{1})-(\beta_{2}x_{1}+\beta_{3}\tilde{x}_{1})\alpha_{3}x_{3}&  \\
\dot{x_{3}}=\alpha_{2}x_{2}+\alpha_{3}\tilde{x}_{2})(\beta_{2}x_{1}+\beta_{3}\tilde{x}_{1})
-(\beta_{2}x_{1}+\beta_{3}\tilde{x}_{1})\alpha_{4}x_{3}-(\beta_{1}x_{2}+\beta_{4}\tilde{x}_{2})^{2}\alpha_{3}x_{3}&  \\
\end{array}%
\right.
\end{equation}
where we considered the following notations:
$$\alpha_{1}=\varepsilon_{0}+\varepsilon_{1}, \,
\alpha_{2}=\varepsilon_{0}+\varepsilon_{2}, \,
    \alpha_{3}=\varepsilon_{1}+\varepsilon_{2},
    \, \alpha_{4}=\varepsilon_{1}+\varepsilon_{3}, \, \alpha_{5}=\varepsilon_{3}+\varepsilon_{2}$$
$$ \beta_{1}=\delta_{0}+\delta_{1}, \, \beta_{2}=\delta_{0}+\delta_{2}, \, \beta_{3}=\delta_{1}+\delta_{3}, \, \beta_{4}=\delta_{2}+\delta_{3}.$$

\begin{re}
\par The analysis of stationary points of the system (\ref{4.14})
is quite difficult, that is why we will present the main results
for fractional Rabinovich differential system.
\end{re}

\section{Fractional Rabinovich differential systems}
\par Generally speaking, there are three mostly used definitions for fractional
derivatives, i.e. Gr$\ddot{u}$nwald-Latnikov fractional
derivatives, Riemann-Liouville fractional derivatives and Caputo's
fractional derivatives, (\cite{Al1},\cite{Di}). Here we discuss
Caputo derivative:
\begin{equation}\label{5.1}
D_{t}^{\alpha}x(t)=I^{m-\alpha}\Big(\frac{d}{dt}\Big)^{m}x(t),
\quad \alpha>0
\end{equation} where $m-1<\alpha\leq m, \, m\geq 1$,
 $\Big(\frac{d}{dt}\Big)^{m}=\frac{d}{dt}\circ...\circ
\frac{d}{dt}$, $I^{\beta}$ is the $\beta^{th}$ order
Riemann-Lioville integral operator, which is expressed in the
following manner:
\begin{equation}\label{5.2}
I^{\beta}x(t)=\frac{1}{\Gamma(\beta)}\int_{0}^{t}(t-s)^{\beta
-1}x(s)ds, \quad \beta>0.
\end{equation}
\par In this paper we consider that $\alpha\in(0,1).$
\par A fractional system of differential equations with
distributed delay in $\mathbb{R}^{3}$ is given by:
\begin{equation}\label{5.3}
D_{t}^{\alpha}x(t)=X(x(t),\tilde{x}(t)), \quad \alpha\in(0,1)
\end{equation} where $x(t)=(x_{1}(t),x_{2}(t),x_{3}(t))\in
\mathbb{R}^{3}.$ \par The matrix associated to the linear part of
the system (\ref{5.3}) in the stationary point $x_{0}$ is given by
the linear fractional differential system:
\begin{equation}\label{5.4}
D_{t}^{\alpha}u(t)=Au(t)+B\tilde{u}(t),
\end{equation}
where $ A=\Big(\frac{\partial X}{\partial x}\Big)\Big|_{x=x_{0}}$
and $ B=\Big(\frac{\partial X}{\partial
\tilde{x}}\Big)\Big|_{x=x_{0}}.$
\par The characteristic equation of (\ref{5.4}) is:
\begin{equation}\label{5.5}
\triangle(\lambda)=det(\lambda^{\alpha}-A-k^{1}(\lambda)B)
\end{equation} where $k^{1}(\lambda)=\int_{0}^{\infty}k(s)e^{-\lambda
s}ds$  and k is given by (\ref{4.4})-(\ref{3.7}).\\*
\begin{pr}\cite{De}
\begin{enumerate}
    \item If all the roots of the characteristic equation
    $\triangle(\lambda)=0$ have negative real parts, then the
    stationary point $x_{0}$ of (\ref{5.5}) is asymptotically
    stable.
    \item If k(s) is Dirac distribution, the characteristic
    equation is given by:
\begin{equation}\label{5.6}
\triangle(\lambda)=det(\lambda^{\alpha}-A-e^{-\lambda \tau}B)=0.
\end{equation} If $\tau=0,$ $\alpha \in (0,1)$ and all the
roots of the equation $det(\lambda I-A-B)=0$ satisfies $\mid
arg(\lambda)\mid>\frac{\alpha \pi}{2},$ then the stationary point
$x_{0}$ is asymptotically stable.
    \item If $\alpha \in (0.5,1)$ and the equation $det(\lambda I-A-e^{-\lambda
    \tau}B)=0$ has no purely imaginary roots for any $\tau>0,$ then
    the stationary point is asymptotically stable.
\end{enumerate}$\Box$
\end{pr}

\par Let us consider a fractional 2-tensor field $P^{\alpha}\in
\mathcal{X}^{\alpha}(\mathbb{R}^{3})\times\mathcal{X}^{\alpha}(\mathbb{R}^{3})$
and $d^{\alpha}f,d^{\alpha}g\in \mathcal{D}(\mathbb{R}^{3}).$ The
bilinear map \\ $[\cdot,
\cdot]^{\alpha}:\mathbf{C}^{\infty}(\mathbb{R}^{3})\times\mathbf{C}^{\infty}(\mathbb{R}^{3})\longrightarrow
\mathbf{C}^{\infty}(\mathbb{R}^{3})$ defined by:
$$[f,g]^{\alpha}=B^{\alpha}(d^{\alpha}f,d^{\alpha}g), \quad f,g\in\mathbf{C}^{\infty}(\mathbb{R}^{3})
$$ is called the fractional Leibniz bracket.

\par If $P^{\alpha}$ is skew-symmetric, we say that $(\mathbb{R}^{3},[\cdot,
\cdot]^{\alpha})$ is a fractional almost Poisson manifold. For
$h\in \mathbf{C}^{\infty}(\mathbb{R}^{3})$ the fractional almost
Poisson dynamical system is given by:
\begin{equation}\label{5.7}
D_{t}^{\alpha}x_{i}(t)=[x_{i}(t),h(t)]^{\alpha}, \quad
[x_{i},h]^{\alpha}=\sum_{i,j=1}^{3}P_{ij}^{\alpha}D_{x_{j}}^{\alpha}.
\end{equation}
\par Let $P^{\alpha}$ be a skew-symmetric fractional 2-tensor field and a symmetric fractional \\* 2-tensor field
$g^{\alpha}$ on $\mathbb{R}^{3}.$ We define the bracket \\
$[\cdot,
\cdot]^{\alpha}:\mathbf{C}^{\infty}(\mathbb{R}^{3})\times\mathbf{C}^{\infty}(\mathbb{R}^{3})\longrightarrow
\mathbf{C}^{\infty}(\mathbb{R}^{3})$ by:
$$[f,g]^{\alpha}=P^{\alpha}(d^{\alpha}f,d^{\alpha}h)+g^{\alpha}(d^{\alpha}f,d^{\alpha}h), \quad
f,h\in\mathbf{C}^{\infty}(\mathbb{R}^{3}).$$

\par The 4-tuple
$(\mathbb{R}^{3},P^{\alpha},g^{\alpha},[\cdot,\cdot]^{\alpha})$ is
called fractional almost metric manifold. The fractional dynamical
system associated to $h\in \mathbf{C}^{\infty}(\mathbb{R}^{3})$
and
\begin{equation}\label{5.8}
D_{t}^{\alpha}x_{i}(t)=[x_{i}(t),h(t)]^{\alpha}, \quad
[x_{i},h]^{\alpha}=\sum_{i,j=1}^{3}P_{ij}^{\alpha}D_{x_{j}}^{\alpha}+\sum_{i,j=1}^{3}g_{ij}^{\alpha}D_{x_{j}}^{\alpha}.
\end{equation}
\begin{pr}
\begin{enumerate}
    \item The fractional dynamical system (\ref{5.7}) is given by:
\begin{equation}\label{5.0}
\left\{%
\begin{array}{ll}
     D_{t}^{\alpha}x_{1}(t)=x_{2}(t)x_{3}(t) &  \\
     D_{t}^{\alpha}x_{2}(t)=-x_{1}(t)x_{3}(t) &  \\
     D_{t}^{\alpha}x_{3}(t)=x_{1}(t)x_{2}(t)&  \\
\end{array}%
\right.
\end{equation}
    \item The fractional dynamical system (\ref{5.8}) is given by:

\begin{equation}\label{5.00}
\left\{%
\begin{array}{ll}
     D_{t}^{\alpha}x_{1}(t)=((\alpha_{1}x_{3}(t)+\alpha_{4}\tilde{x}_{3}(t))(\beta_{1}x_{2}(t)+\beta_{4}\tilde{x}_{2}(t))\\*
+(\beta_{2}x_{1}(t)+\beta_{3}\tilde{x}_{1}(t))(\beta_{1}x_{2}(t)+\beta_{4}\tilde{x}_{2}(t))\alpha_{2}x_{3}(t) &  \\
     D_{t}^{\alpha}x_{2}(t)=-(\alpha_{1}x_{3}(t)+\alpha_{4}\tilde{x}_{3}(t))(\beta_{2}x_{1}(t)+\beta_{3}\tilde{x}_{1}(t))-(\beta_{2}x_{1}(t)+\beta_{3}\tilde{x}_{1}(t))\alpha_{3}x_{3}(t) &  \\
     D_{t}^{\alpha}x_{3}(t)=(\alpha_{2}x_{2}(t)+\alpha_{3}\tilde{x}_{2}(t))(\beta_{2}x_{1}(t)+\beta_{3}\tilde{x}_{1}(t))\\*
-(\beta_{2}x_{1}(t)+\beta_{3}\tilde{x}_{1}(t))\alpha_{4}x_{3}(t)-(\beta_{1}x_{2}(t)+\beta_{4}\tilde{x}_{2}(t))^{2}\alpha_{3}x_{3}(t).&  \\
\end{array}%
\right.
\end{equation}

    \item The fractional dynamical systems (\ref{5.0}) and
    (\ref{5.00}) have the stationary points $e_{1}^{m}(m,0,0)$, $e_{2}^{m}(0,m,0)$ and
    $e_{3}^{m}(0,0,m),$  $m \in \mathbb{R}.$
    \item The characteristic equations for (\ref{5.0}) are given by:\begin{itemize}
        \item $e_{1}^{m}(m,0,0):$ $\lambda^{\alpha}(-\lambda^{2\alpha}+m^{2}(m+1))=0;$
        \item $e_{2}^{m}(0,m,0):$ $\lambda^{\alpha}(\lambda^{2\alpha}+m^{2}\lambda^{2}-m^{2})=0;$
        \item $e_{3}^{m}(0,0,m):$ $\lambda^{\alpha}(\lambda^{2\alpha}+m^{2})=0.$
    \end{itemize}
    \item The characteristic equations for (\ref{5.00}) are given by:\begin{itemize}
        \item $e_{1}^{m}(m,0,0):$ $\lambda^{\alpha}(\lambda^{2\alpha}+a\lambda^{\alpha}+be^{-\lambda^{\alpha} \tau}+c)=0, \quad a,b,c \in \mathbb{R};$
        \item $e_{2}^{m}(0,m,0):$ $\lambda^{\alpha}(\lambda^{2\alpha}+a_{1}\lambda^{\alpha}+b_{1} e^{-\lambda^{\alpha} \tau}+c_{1}e^{-2\lambda^{\alpha}\tau})=0, \quad a_{1},b_{1},c_{1}\in \mathbb{R};$
        \item $e_{3}^{m}(0,0,m):$ $\lambda^{\alpha}(\lambda^{2\alpha}+a_{2}\lambda^{\alpha}+b_{2} e^{-\lambda^{\alpha} \tau}+c_{2}e^{-2\lambda^{\alpha}\tau}+d_{2})=0, \quad a_{2},b_{2},c_{2},d_{2}\in \mathbb{R}.$
    \end{itemize}
\end{enumerate}
\end{pr}

 \par In the second section, "The analysis of classical Rabinovich differential
equations", we worked with a single-step method, Runge-Kutta, that
means that we used only the information regarding the previous
point for computing the successive point. For our illustrations we
develop the Adams-Bashforth-Moulton predictor-corrector method.
\par We can integrate numerically the set of differential
equations (\ref{5.0}) with the Adams-Bashforth-Moulton method. To
do so, we consider the relations we used for our simulation part,
in Moulton method:
\begin{eqnarray*}
x_{1}(j+1)&=&x_{1}(0)+\frac{1}{\Gamma(\alpha)}(\sum_{k=0}^{j}a(k,
j+1)x_{2}(k)x_{3}(k)+a(j+1,j+1)x_{2p}(j+1)x_{3p}(j+1))\\
x_{1p}(j+1)&=&x_{1}(0)+\frac{1}{\Gamma(\alpha)}(\sum_{k=0}^{j}b(k,
j+1)x_{2}(k)x_{3}(k))\\
x_{2}(j+1)&=&x_{2}(0)+\frac{1}{\Gamma(\alpha)}(\sum_{k=0}^{j}-a(k,
j+1)x_{1}(k)x_{3}(k)-a(j+1,j+1)x_{1p}(j+1)x_{3p}(j+1))\\
x_{2p}(j+1)&=&x_{2}(0)-\frac{1}{\Gamma(\alpha)}(\sum_{k=0}^{j}b(k,
j+1)x_{1}(k)x_{3}(k))\\
x_{3}(j+1)&=&x_{3}(0)+\frac{1}{\Gamma(\alpha)}(\sum_{k=0}^{j}a(k,
j+1)x_{1}(k)x_{2}(k)+a(j+1,j+1)x_{1p}(j+1)x_{2p}(j+1))\\
x_{3p}(j+1)&=&x_{3}(0)+\frac{1}{\Gamma(\alpha)}(\sum_{k=0}^{j}b(k,
j+1)x_{1}(k)x_{2}(k)),\\
\end{eqnarray*}
where $$i=0,1,..., n \quad and \quad j=0,1,...,m$$ and
\begin{eqnarray*}
b(i,j+1)&=&h^{\alpha}\frac{(j-i+1)^{\alpha}-(j-i)^{\alpha}}{\alpha}\\
a(i+1,j+1)&=&h^{\alpha}\frac{(j-i+1)^{\alpha
+1}+(j-i-1)^{\alpha+1}-2(j-i)^{\alpha
+1}}{\alpha(\alpha +1)}\\
a(0,j+1)&=&h^{\alpha}\frac{j^{\alpha +1}-(j-\alpha)(j+1)^{\alpha}}{\alpha(\alpha +1)}\\
a(j+1,j+1)&=&\frac{h^{\alpha}}{\alpha(\alpha+1)}.
\end{eqnarray*}

\par We consider a graphic representation of Moulton method, for
our system (\ref{1.1}), for the following two cases:
\par (1) $x_{1}(0)=0.001$, $x_{2}(0)=0.001,$ $x_{3}(0)=6$, $\alpha=0.8$ (Figure 1);
\par (2) $x_{1}(0)=0.001$, $x_{2}(0)=0.001$, $x_{3}(0)=6$, $\alpha=1$ (Figure 2)
$$ $$
$$ $$
$$ $$
$$ $$
$$ $$
$$ $$
$$ $$
\begin{figure}
 \includegraphics[height=3in]{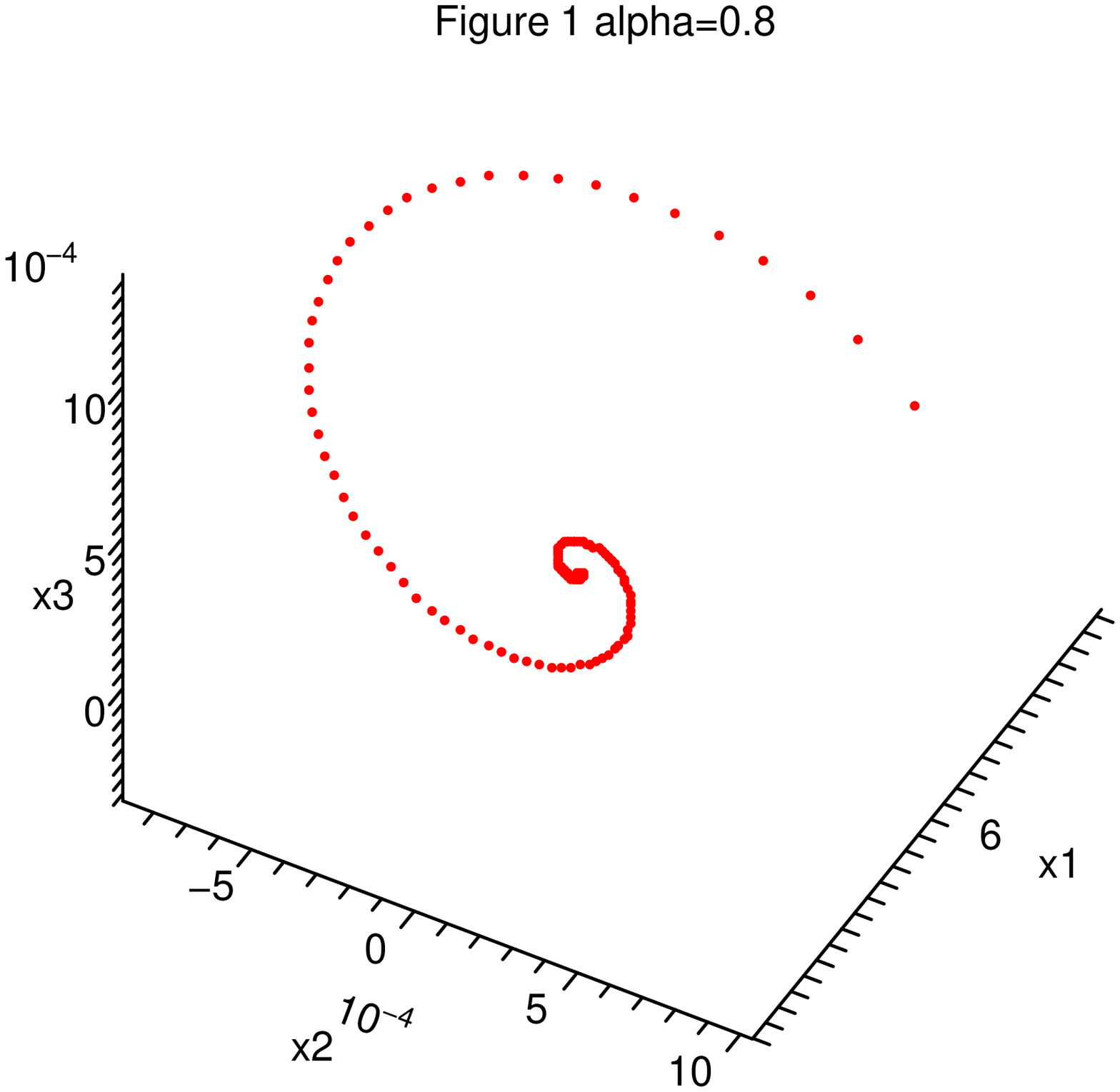}
 \includegraphics[height=3in]{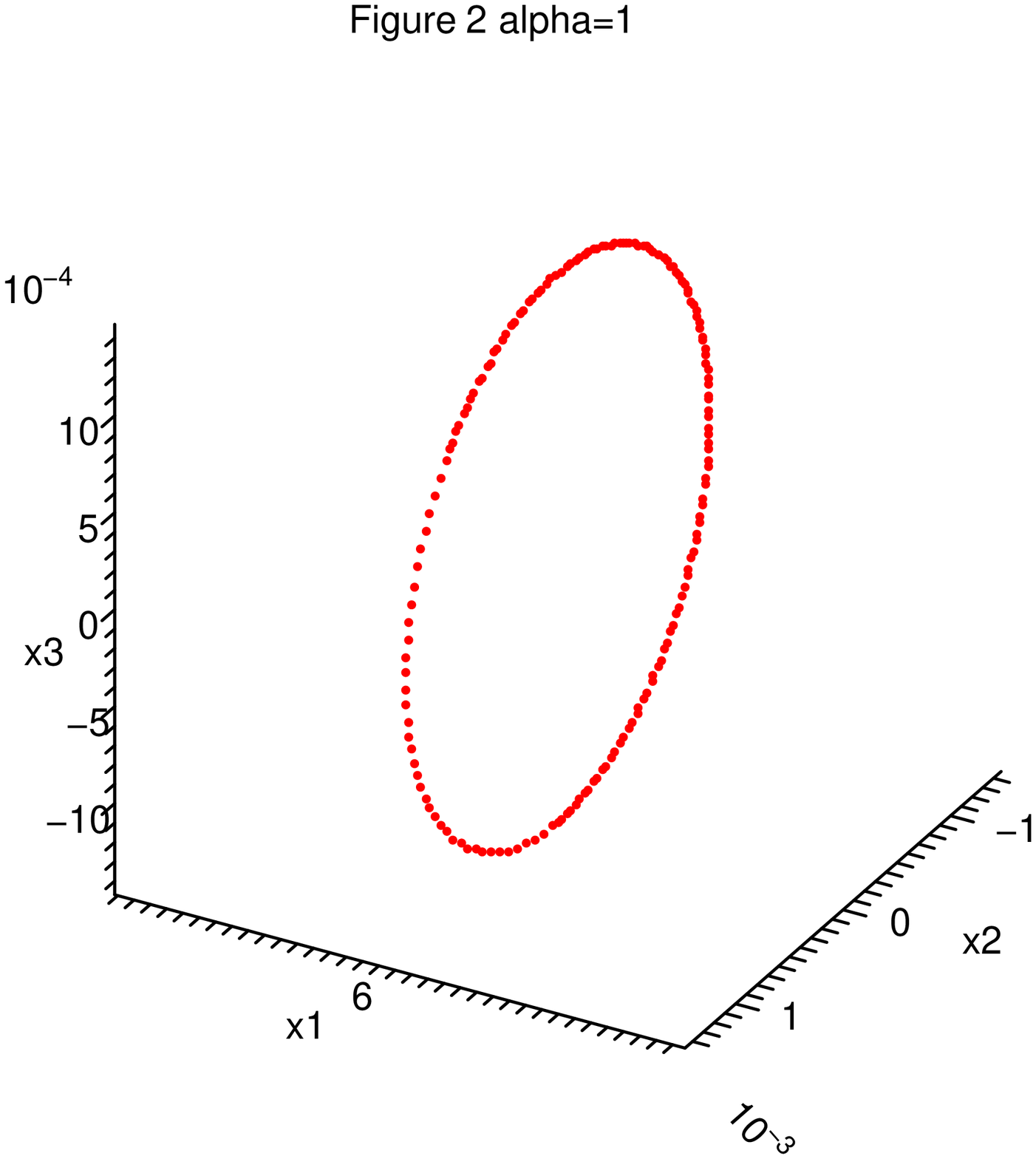}
  \end{figure}

\newpage

\section{Conclusions}
\par Until now we have an approach and also some solutions for
metriplectic manifolds of the first and second kind, and also
differential systems with distributed delay(for the first case
presented here), and Rabinovich fractional differential system,
with Dirac distribution$(\tilde{x}(t)=x(t-\tau))$.What we want to
continue is to apply all other distributions for our three cases.
\\*

\underline{Acknowledgements }: The authors were partially
supported by the Grant CNCSIS 95GR 2007/2008.

\center
 Seminar of Geometry-Topology\\
 West University of Timi\c{s}oara\\
 B-dul V.P\^{a}rvan no 4,\\
 300223 Timi\c{s}oara, Romania\\
 email: chisoana@yahoo.com\\
 email: puta@math.uvt.ro


\begin{thebibliography}{99}
\bibitem{Al1} I.D. Albu, M.Neantu, D.Opris, \emph{The geometry of fractional of osculator bundle of higher order and
applications}, Conference of Differencial Geometry: Lagrange and
Hamiltoniah Spaces, Sepember, 3-8, 2007, Iasi.

\bibitem{Al2} I.D.Albu, M.Neamtu, D.Opris, \emph{Dissipative mechanical systems with
delay}, Tensor N.S. vol 67(2006), 1-27.

\bibitem{Al3} I.D. Albu, D.Opris, \emph{Leibniz dynamics with time
delay}, arXive.math/0508225, 15p.

\bibitem{Be} B.H. Bermejo and V. Fairen, \emph{Simple evaluation of Casimir invariants in finite dimensional  Poisson
systems}, Phys. Lett. A 241 (1998), 148-154.

\bibitem{Bi1} P.Birtea, M. Puta and R.M. Tudoran, \emph{Periodic orbits in the case of a zero eigenvalue} (to appear in C.R. Acad. Sci.
Paris).

\bibitem{Bi2} P. Birtea, M. Boleantu, M. Puta, R.M. Tudoran,
\emph{Asymtotic Stability for a Class of Metriplectic Systems},
arXiv:071.3012v1 [math-ph] 16 Oct 2007.

\bibitem{Ch} O. Chis, M. Puta, \emph{Geometrical and dynamical aspects in the theory of Rabinovich system},
 The eigth international workshop on differential geometry and its applications, August 19-25, 2007, Cluj-Napoca,
 Romania.

\bibitem{De} W. Deng, C. Li, J. L$\ddot{u}$, \emph{Stability analysis of linear fractional differential system with multiple time
 delay}, Nonlinear syn(2007)48:409-416.

\bibitem{Di} K. Diethem, \emph{Fractional Differential Equations, Theory and Numerical
 Threatment}, Braunschweig, 2003.

\bibitem{Fi} D. Fish, \emph{Dissipative perturbation of 3D Hamiltonian
systems}, arXive: math.ph/0506047. v1, 2005, 12 pg.

\bibitem{Or} J.P. Ortega, V. Planas-Bielsa, \emph{Dynamics on Leibnitz
manifolds}, arXive: math DS/0309263, 2503.

\bibitem{Pi} A.S. Pikovsky, M.I. Rabinovich, Math. Phys. Rev. 2,
165, (1981).

\bibitem{Pu1} M.Puta, P. Birtea and R.M.Tudoran, \emph{Poisson manifolds and Bermejo-Fairen construction of
Casimirs}, Tensor N.S. vol 66 (2005) 59-70.

\bibitem{Pu2} M.Puta,
\emph{Hamiltonian systems and geometric quantisation}, Mathematics
and Applications vol 260, Kluwer Academic Publishers, 1993.\\*
\end{thebibliography}
\end{document}